\documentclass[11pt]{article}

\usepackage{graphicx}
\usepackage{algorithmic}
\usepackage[lined,boxed,commentsnumbered]{algorithm2e}
\usepackage{amssymb}
\newcommand{\Bbbr}{\mathbb{R}}
\newtheorem{definition}{Definition}[section]

\newtheorem{Theorem}[definition]{Theorem}

\begin{document}
\title{How to get high order without loosing efficiency for the resolution of systems of nonlinear equations: A short review of Shamanskii's $m$ method
}

\author{F. Calabr\`o \footnote{DIEI, Universit\`a degli Studi di Cassino e del Lazio Meridionale, Cassino (FR), Italy. E-mail: calabro@unicas.it} \& A. Polsinelli\footnote{a\_polsi@hotmail.it} }

\maketitle

\begin{abstract}
We present relations between some recently proposed methods for the solution of a nonlinear system of equations. In particular, we review the Shamanskii's $m$-method \cite{shamanskii1967modification}, that is an iterative method derived from Newton's method that converge with order $m+1$. We discuss efficient implementation of this method via matrix factorization and some relevant properties.
\\
We believe that recent developments in the research of solutions of systems of equations did not take sufficiently into account this method. The hope, with this paper, is to encourage the entire community to remember this simple method and use it for comparison when new methods are introduced.
\\
\textit{This work is dedicated to Prof. Elvira Russo: a very special teacher.}

\end{abstract}

\section{Introduction}
The design of iterative algorithms for solving systems of nonlinear equations is a numerical analysis topic that has never known years with little interest. A result is that a researcher approaching such field and trying to get a complete review of the main developments has the feeling that this is close to impossible. When interested in a method for the search of a solution of system of nonlinear equations, one has to fix clearly the hypotheses that are known in the specific problem and work with these as guidelines. Our study has been motivated by the construction of quadrature rules with ad-hoc exactness, as considered in \cite{bremer2010nonlinear,ma1996generalized} to be used for some new applications related to isogeometric analysis, see \cite{renePRE,Johannessen:2016,manni2011generalized,aimi2016isogemetric}. In our application the problem is squared and exact Jacobian is available; moreover the problem is of medium size (not more than $50$ equations) and a good starting point is computable. These are the main features that we consider, for details we refer to the cited papers. \\
In the original papers \cite{Hughes:2010,ma1996generalized}, a Newton procedure is used for the resolution of the nonlinear equations that fix the quadrature rule, but in many cases the procedure fails to converge, see \cite{Hughes:2010,Johannessen:2016}. While searching for a good strategy to solve our problem, we have encountered many methods that modify Newton's iteration and converge with higher order, see\footnote{The list is far from complete: also many methods have been introduced for the resolution of the univariate problems but can be easily extended to the multivariate case.} e.g.  \cite{darvishi2007third,frontini2003some,kou2006modification,mcdougall2014simple,homeier2004modified}. All these methods are referred to as Modified-Newton: their main feature is that the Jacobian can be evaluated so that there is no need for an approximation as for Quasi-Newton methods. When counting FLOPs, the most expensive part of Modified-Newton methods is where the inversion of the Jacobian is performed: in all the cited papers the number of matrix inversion is raised, so that efficiency is reduced. All modifications of Newton's method consider more inversion of the Jacobian matrix. Nevertheless, many of these methods have remarkable interesting features, see for example the one derived by quadrature rules applied to the integral form of the reminder term in Taylor's formula in \cite{weerakoon2000variant,frontini2004third,gracca2015root}.\\
%, so that the order of convergence is   All this modifications have some interesting feature and some are really difficult to distinguish each other. 
We  focus in this paper on the family of Shamanskii's $m$ methods. The interesting feature of these methods is that at each iteration only one matrix inversion is required, while theoretical convergence rate can be raised thanks to internal iterations of lower computational cost. 
\\
Moerover, with a good choice of the paramenter $m$ the convergence domain  \cite{ezquerro2013local} of Shamanskii's $m$ method seems to be greater of the one of Newton's method.
\\ 
While previously considered a standard procedure (see \cite{brent1973some} and the book \cite{kelley2003solving}) recently Shamanskii's methods seem somehow not prefered, with some noticeable\footnote{We notice here that the third order method in \cite{darvishi2007third} is exactly Shamanskii's method with $m=2$.} exception \cite{kchouk2014new}.  
\\
In this paper we test Shamanskii's method on some examples taken from  \cite{cordero2012increasing,brent1973some}. In all test with non-singular Jacobian the method converges in few iterations. The tests confirm -at least up to order 4- that the order is increased when internal iterations are raised. Moreover, because of the limited number of inversion of matrices considered, we notice that the method is efficient.
\\
The paper is organized as follows. In Section 2 we review Shamanskii's $m$ methods, give a pseudo-code and some details on the convergence properties; in Section 3 we perform some tests and give some conclusions.

\section{Shamanskii's $m$ methods}
\begin{algorithm} \LinesNumbered \SetKwInOut{Input}{Input}
\Input{Starting point $x_0$, function $F$, Jacobian $J$ }
Inizializations: $x^{(0)}= \tilde{x}^{(0)}= x_0$, $it_{tot}=it_{inv}=0$, $rhs=F(\tilde{x}^{(0)})$, $res^{(0)}=norm(rhs)$.\;
\While{Convergence}{ 
$[L,U]=lu(J(x^{(it_{inv})} ))$ \,,\ $it_{inv}=it_{inv}+1$ \;
\For{ $k=1:m$ }{ 
$\tilde{x}^{(it_{tot}+1)} = \tilde{x}^{(it_{tot})} -U\backslash (L\backslash rhs) $ \;
$rhs= F(\tilde{x}^{(it_{tot}+1)}) $ \,,\ $it_{tot}=it_{tot}+1$ 
}
$x^{(it_{inv})}= \tilde{x}^{(it_{tot})}$\,,\ $res^{(it_{inv})}=norm(rhs)$
}
\caption{Shamanskii's $m$-method}\label{algo_iniz}
\end{algorithm}
Shamanskii's $m$ methods are designed for the resolution of a system of nonlinear equations $F(x)=0$ where $F: x\in \Bbbr^n \to F(x) \in \Bbbr^n$ is a functional and $J:  x\in \Bbbr^n \to J(x) \in \Bbbr^{n\times n}$ is it's Jacobian, namely $J_{ij}:= \frac{\partial F_i}{\partial x_j}$. The method computes a sequence of points $x^{(it)}$ by applying a step of Newton's method and $m-1$ internal iterations with frozen Jacobian, namely chord's iterations. We summarize in Algorithm \ref{algo_iniz} the flow-chart, notice that one matrix inversion is needed at each internals iteration so that factorization is helpful in this case. In the next theorem we summarize the properties of convergence, see for example  \cite[Section 2.4]{kelley2003solving}. For similar results with inexact Jacobian -in the case of Quasi-Newton's method- we refer to \cite{brent1973some}.

\begin{Theorem}\label{Th1}
Let $F$ be differentiable and call $J$ it's Jacobian. Suppose that:
\begin{itemize}
\item $x^*\in \Bbbr^n$ exists such that $F(x^*)=0$;
\item $\|J(x)-J(y) \| \le \gamma \|x-y\| $ for some $\gamma$;
\item $J(x^*)$ is positive definite.
\end{itemize}  
Then, for all $m\in\mathbb{N}$ there are $K > 0 $ and $\delta > 0$ such that if $\|x_0 - x^*\|\le \delta$ the Shamanskii's iterates $x^{(j)}$ converge to $x^*$ with order $m+1$:
\begin{equation}
\| x^{(j+1)}-x^*\| \le K \| x^{(j)}-x^*\|^{m+1} .
\end{equation}
\end{Theorem}

Theorem \ref{Th1} ensures  that within the family of Shamanskii's method it is possible to achieve a desired convergence order by simply increasing $m$, the number of internal iterations. Obviously the limit case ($m\to \infty$) would give an ''infinite order'' method while it simply coincides with the chord method of order 1. This ambiguity is due to the conventional assumption that counts as iteration only the iterations embedding a matrix inversion. 
In next section we report the results of our testing accounting both for outer and internal iterations, since  most of the computational costs are incurred within the outer ones, mainly for large problems (i.e. $n>>m$). The choice of the parameter $m$ to be used\footnote{Some discussion on this is available in \cite{brent1973some} in the case of Quasi-Newton implementation.} strongly depends on the problem, and in authors' opinion some effort should be spent for an implementation that adaptively interrupts the internal iterations.

\section{Numerical results}
To test order of convergence of our procedure, we use the usual estimate of convergence order (\cite{cordero2012increasing,brent1973some}) :
\begin{equation}\label{rho}
\rho= \frac{\log\left( \frac{ \| x^{(j)} - x^{(j-1)}\| }{ \| x^{(j-1)} - x^{(j-2)}\| } \right)}  { \log\left( \frac{ \| x^{(j-1)} - x^{(j-2)}\| }{ \| x^{(j-2)} - x^{(j-3)}\| }\right)} 
\end{equation}

\begin{table}
{\small
\begin{tabular}{c|rrrrrrrr}
 & $m=1$ &2 \rule{0.6cm}{0cm} &3 \rule{0.6cm}{0cm} &4 \rule{0.6cm}{0cm}  \\
\hline
(a) &  5 (5) 2.0044 &    3 (6) 2.9649&    3 (9) 3.9146&   2(8) NA   \\
(b) &  6 (6) 1.9946 &    4 (8)  2.9593&    3 (9) 3.3058&   3(12) 4.2046 \\
(c) &  5 (5) 1.9127 &    3 (6)  2.7689&    3 (9) 3.7390&   3(12)  4.7116 \\
(d) &  6 (6) 1.9968 &    4 (8) 2.9594 &    3 (9)  3.3493 &   3(12) 4.2464 \\
(e) &  7 (7) 2.0002 &    5 (10) 2.9754 &    5 (15) 3.7121 &   6(24) 4.6198 \\
\end{tabular}
}
\caption{Convergence test. For each test problem in (\ref{examples})  we report the number of iterations $it_{inv} $ needed by Shamanskii's $m$ method for convergence (in parentheses the total number of iterations $it_{tot}$) and the computed order of convergence $\rho$ that follows from equation (\ref{rho}).}
\label{iter}
\end{table}

In table \ref{iter} we report the computed $\rho$ and the number of iteration made to meet the tolerance on the residual, fixed to $10*eps$. We have considered the following functionals and starting points taken from \cite{cordero2012increasing,brent1973some}:
\begin{equation}\label{examples}
\begin{array}{c|l|c}
&F(x) & x_0 \\[5pt] \hline & & \\
(a) &[x_1^2-4x_2+x_2^2 ; 2x_1-x_2^2 - 2] & [1;0.1] \\
(b) & [ x_1^2+x_2^2-1 ; x_1^2-x_2^2 + 0.5 ] & [1;1] \\
(c) & [ \cos(x_2)-\cos(x_1) ; x_3^{x_1}-1/x_2; exp(x_1)-x_3^2] & [1;1;2]\\
(d) & [x_i x_{i+1}-1 _{i=1,\dots 30} ; x_{31}x_1 -1 ] &-2 ones(31,1)\\
(e) & [x_1^2+x_2^2-2 ; exp(x_1-1)+ x_2^2-2] & [2;0.5]
\end{array}
\end{equation}
The convergence test confirm that the Shamanskii's $m$ method increases the order of convergence when increasing $m$; correspondingly the number of iterations with matrix inversion are lowed and the overall computational efficiency is raised. 
Our conclusion is that Shamanskii's method is both efficient and robust, so that it is still the Modified-Newton method that should be used -at least as a comparison- when higher order is required.
\section*{Acknowledgements}
Francesco Calabr\`o was partially supported by INdAM, through GNCS research projects. This support is gratefully acknowledged.

\bibliographystyle{elsarticle-num}
\bibliography{NW}

%\section*{Figures and tables}

\end{document}